\documentclass[12pt]{amsart}

\def\checkbox{\leavevmode\vbox to 9pt{\hrule \vss
	\hbox to 9pt{\vrule height 9pt \hfil\vrule height 9pt}\vss
	\hrule}\ }

\newcommand{\Q}{\mathbb Q}

\newcommand{\N}{\mathbb N}

\renewcommand{\epsilon}{\varepsilon}
\renewcommand{\phi}{\varphi}
\newtheorem{Lemma}{Lemma}[section]
\newtheorem{Theorem}{Theorem}[section]
\newtheorem{Proposition}{Proposition}[section]

\newtheorem{Corollary}[Proposition]{Corollary}

\newtheorem{Conjecture}{Conjecture}[section]

\begin{document}
\begin{abstract} We prove that all Fano threefolds with log-terminal singularities of given index belong to finitely many families. This result was previously obtained by the author in the case of unipolar Fano varieties.
\end{abstract}

\author[]{Alexandr Borisov}
\title[A. Borisov, Fano threefolds]{Boundedness of Fano threefolds with log-terminal singularities of given index}
\address{Department of Mathematics, Penn State University, University Park, PA 16802, USA}
\email{borisov@math.psu.edu}
\maketitle

\centerline{borisov@math.psu.edu}


\section{Introduction}                                          

This paper is concerned with Fano varieties with log-terminal singularities.
Fano varieties and log-terminal singularities play very important role in
modern birational algebraic geometry (cf. \cite{CKM}, \cite{IskovskikhProkhorov}, \cite{KMM}).

The main goal of this paper is to prove the following theorem.

\begin{Theorem} For any $n\in \N$ the family of Fano threefolds with log-terminal singularities of index $n$ is bounded.
\end{Theorem}

This theorem is the three-dimensional case of the following 
conjecture originally proposed by V. Batyrev.
It was previously obtained by the author in \cite{Fano} under
the additional assumptions that $X$ is unipolar, i.e. it is $\Q-$factorial and the Picard 
number $\rho (X)$ is equal to $1.$

\begin{Conjecture}(Batyrev)
For any two fixed natural numbers $k,$ $n$ there are just finitely many 
families of Fano varieties
of dimension $k$ with log-terminal singularities of index $n$. 
\end{Conjecture}

Theorem 1.1. has several noteworthy corollaries.

\begin{Corollary} The family of Fano threefolds with rational Gorenstein singularities is bounded.
\end{Corollary}

\begin{Corollary} Suppose $X$ is a Fano threefold with log-terminal singularities. Then the algebraic fundamental group of its smooth locus is finite.
\end{Corollary}

Also, combined with the recent work of Koll\'ar, Miyaoka, Mori and Takagi it implies boundedness of Fano threefolds with canonical singularities (cf. \cite{KMMT}).

Through the above application and by itself, Theorem 1.1 is a step towards proving the following stronger conjecture, proposed
independently by the author (cf. \cite{Fano}) and V. Alexeev
(cf. \cite{Alexeev2}).

\begin{Conjecture}(Borisov-Alexeev) 
For any fixed $n\in \N,$ $\epsilon >0$ there are just finitely
many families of $n-$dimensional Fano varieties with 
$\epsilon-$logterminal singularities.
\end{Conjecture}

This conjecture is important for the so-called log Sarkisov Program (cf. \cite{BrunoMatsuki}).
Some partial results towards it are due to Alexeev and Nikulin
(in dimension 2), the author and L. Borisov (in the toric case), and
Y. Kawamata (cf. \cite{Kawamata}).  The boundedness of smooth Fano varieties
in any given dimension was proven by Koll\'ar, Miyaoka, and Mori (cf. 
\cite{KoMiMo1}). Their method (the ``rational curve surgery") was the basis
of the method of \cite{Fano}. Recently, Koll\'ar, Miyaoka, Mori and Takagi proved boundedness of three-dimensional Fano varieties with canonical singularities (using the result of this paper).  Also recently, an entirely new approach to the boundedness of Fano varieties was proposed by Ziv Ran (cf. \cite{Ran}, \cite{RanClemens}). It is based on studying sheaves of high order differential operators on the plurianticanonical bundles, and has some potential to compete with the rational curve surgery.

Also, H. Tsuji recently announced the proof of Batyrev conjecture in any dimension, under extra assumptions that $X$  is unipolar (cf. \cite{Tsuji}). His methods are analytic, and are yet to be independently verified. 

In this paper we employ the techniques of \cite{Fano} together with some
new ideas. These new ideas include the following.

\begin{enumerate}
\item Using Alexeev Minimal Model program with suitable boundary to find
horizontal extremal contractions.

\item Koll\'ar's effective Base Point Freeness theorem.

\item Using Kawamata's result on the length of extremal curves with suitable 
boundary to avoid gluing curves in some cases.
\end{enumerate}

The paper is organized as follows. In section 2 we consider the problem
in dimension 2. Of course, Batyrev conjecture in dimension 2 was first
proven by Nikulin  in \cite{Nikulin1}. But the rational curve surgery provides a very easy proof. Also in section 2 we discuss some reformulations of
the problem that follow from Koll\'ar's Effective Base Point Freeness theorem.

In section 3 we prove that, except for a bounded family, all Fano threefolds 
with log-terminal singularities of given index $n$ 
contain a covering family of rational curves $\{l\}$ such that $l\cdot (-K_X)$
is bounded in terms of $n,$ and the rationally connected fibration associated to $\{l\}$ has fibers of
dimension $2$.

In section 4 we complete the proof of Batyrev conjecture in dimension three.

{\bf Notations.} We utilize the notations of \cite{Fano}.
That is, we will usually identify curves on different birationally 
equivalent varieties if they coincide in their general points. The 
identified curves will be usually denoted by the same symbol. The same 
convention will be used for the two-dimensional subvarieties. If it is 
necessary to point out that, say, a prime divisor $S$ is considered on 
a variety $X$ it will be denoted by $S_X$. Another convention is that 
$\{l\}$ will denote the family of curves with a general element $l$ and 
$\{H\}$ will denote the LINEAR system of Weil divisors with a general 
element $H$. Whenever we have a family of curves, a general element is 
reduced and irreducible unless the opposite is explicitly specified. 
It will be clear in every particular case why these conventions agree 
with each other.

Additionally, by a constant $c$ or $c(n)$  we will mean some positive constant
that only depends on $n$.  The value of $c$ can be different in different 
parts of the proof.

{\bf Acknowledgments.} The author is indebted to V. Iskovskikh who introduced him to algebraic geometry and Fano varieties in particular. The author would like to thank V. Alexeev, V. Batyrev, Y. Kawamata, K. Matsuki, S. Mori, M. Reid, and V. Shokurov for their interest in author's birational geometry research over the years. The author is thankful to J. Koll\'ar, J. McKernan and the referee for helpful advice. This research was conducted while the author was working at Washington University, St. Louis. The author would like to thank the algebraic geometers there (N. M. Kumar, V. Masek, D. Wright) for their support and hospitality.


\section{Two-dimensional case and some easy equivalences}
We begin with the following theorem, which is the two-dimensional case of the conjecture of Batyrev, mentioned above.

\begin{Theorem} For any $n\in \N$ the family of log del Pezzo surfaces of index $n$ is bounded.
\end{Theorem}

{\bf Proof.} Suppose $X$ is a log del Pezzo surface of index $n$. By the result of Miyaoka and Mori (cf. \cite{MiyaokaMori}), there
exists a covering family $\{l\}$ of rational curves on $X$ such that
$l\cdot (-K_X) \le 4.$ If a general $l$ passes through some singularity of $X$ then by Lemma 2.2 of \cite{Fano}, $(-K_X)^2$ is bounded. This is because in
dimension $2$ for any log-terminal singularity $x_0\in X$
$mult(x_0) \le 2 \cdot index (x_0).$

If a general $l$ is contained in $X\setminus Sing(X),$ consider the 
minimal resolution (or terminal modification) $\pi :Y\rightarrow X.$ 
There are two cases.

1) Two general points of $Y$ can be connected by a chain of at most two curves from $\{l\}$. 

2) Two general curves from $\{l\}$ do not intersect each other. In this case
$l^2=0,$ $-K_X\cdot l=-K_Y\cdot l=2.$

In the first case $(-K_X)^2\le 2!(2\cdot 4)^2 =128.$ (cf., e.g., \cite{KoMiMo2}).

In the second case, because all log del Pezzo surfaces are rational, there 
exists a rational curve $L$ on $Y$ such that $L\cdot l >0.$ Applying Mori's 
bend-and-break procedure to it, we can find a new rational curve $L^{'},$ 
such that
$L^{'} \cdot l >0$ and $L^{'}$ does not admit a non-trivial deformation with two fixed points. This means, in particular, that $L^{'}\cdot (-K_Y) \le 3.$ Also, $L^{'}\cdot K_Y \le L^{'}\cdot K_X.$ So $L^{'}\cdot K_Y \le -1.$  By one of the gluing lemmas of Koll\'ar-Miyaoka-Mori (cf. \cite{KoMiMo1}) we can glue together $L^{'}$ and at most 2 copies of $l$ to obtain a new family of rational curves $l^{'},$ which connects two general points of $Y$. Therefore $(-K_X)^2 \le (3+2\cdot 2)^2=49.$

So in all cases $(-K_X)^2$ is bounded. The family of such log del Pezzo surfaces is bounded by the results of Koll\'ar (cf. \cite{Kollar1}). \hfill \checkbox

{\bf Remark.}  Instead of applying the gluing lemma above, one can use a simple Riemann-Roch argument. We chose the above proof only for its similarity with our proof of the three-dimensional case.

The following theorem is  a corollary  of the results of Koll\'ar.

\begin{Theorem} Suppose $X$ is  a Fano variety of dimension $k$ with 
log-terminal singularities of index $n.$ Then we have the following.
\begin{enumerate}

\item For some natural number $N$ that only depends on $k$, the Cartier
divisor $-N\cdot n \cdot K_X$ is very ample. For $k=3$ one can choose
$N=4320.$

\item For every $k$  the following
two statements are equivalent. 

a) The family of such varieties is bounded.

b) The ``degree" $(-K_X)^k$ is bounded.
\end{enumerate}

\end{Theorem}

{\bf Proof.} The first assertion is a direct corollary of the Effective Base
Point Freeness theorem of Koll\'ar (cf. \cite{Kollar1},  1.1,  1.2). For the
second assertion we only need to prove that (b) implies (a). If (b) is
satisfied then by the result of Koll\'ar and Matsusaka (cf. \cite{KollarMatsusaka}) all coefficients of the Hilbert
polynomial of $(-nK_X)$ are bounded. Together with the first assertion this
implies the result. \hfill \checkbox

The following theorem shows that instead of bounding $(-K_X)^k$ one can also 
bound the dimension of the space
of global sections of big enough multiple of $(-nK_X).$

\begin{Theorem}  Suppose $X$, $n,$ $k$ are as above. Then there is a constant
$A$ that only depends on $k$ with the following property.

For all $l\ge 2k,$ $(-nK_X)^k\le A\cdot h^0(l\cdot(-nK_X))$
\end{Theorem}

{\bf Proof.} By the Kawamata-Viehweg vanishing theorem (cf. \cite{KMM}), for all $i\ge 0$
$h^0(-inK_X) =\chi (-inK_X)$.

So $f(i)=h^0(-inK_X)$ is a polynomial of degree $k$ for $i\ge 0.$ It obviously 
has  the following five properties.
\begin{enumerate}
\item $f(0)=1$

\item$f(l)=h^0(l\cdot(-nK_X))$

\item $f(i) \ge 0$ for all $i \ge 0$

\item If for some $i\le l$, $f(l-i) >0,$ then $f(i) \le f(l)$

\item The highest degree term of $f$ is $\frac{(-nK_X)^k}{k!}x^k$
\end{enumerate}

Because $f$ is not identically zero, it has at most $k$ integer roots.
So for   any fixed $l\ge 2k$ one can  find a sequence of $k$ integers
$1\le i_1<i_2<...<i_k\le 2k$ such that for all $j$, $f(l-i_j) >0.$ By the
property (4) this implies that $f(i_j)\le f(l)$. Also $f(0)=1\le f(l),$ because
otherwise $f$ would vanish at points $i_j$ and $l$. The only way it could
happen without getting $(k+1)$ zeroes is if $i_k=2k=l$. But in this case
$f(k)=0$ and $f$ vanishes for at least half of the $(2k-2)$ numbers 
$1,2,...,(k-1),(k+1),...,(2k-2), (2k-1).$  So we get $(k+1)$ zeroes anyway.

There are just finitely many possibilities for the sequence $\{i_j\}.$ For
each of them, $f(x)$ is determined by its values at $0$ and $i_j,$
$j=1,2, ...k.$ And we can use the above inequalities and Lagrange interpolation formula at the above points to estimate the highest
coefficient of $f(x)$ above by some multiple of $f(l)$.
Combining these estimates with the property (5) proves the theorem. \hfill \checkbox

{\bf Remark.}  The constant $2k$ in the above theorem can probably be
improved. See, in particular, \cite{Fano}, Lemma 2.1. It would be more 
important, however, to drop the dependence on $n.$ The problem is, of 
course, that this would require some Riemann-Roch theorem for Weil divisors, 
and $h^0(-iK_X)$ is in general not a polynomial for $i\ge 0.$

By the above theorems, to prove the Batyrev conjecture in dimension 3 it is
enough to obtain a bound on the self-intersection or the dimension of
the space of global sections of the very ample divisor $H=-4320nK_X.$ 
This is going to be our goal. Many of our constructions will be carried
out under the assumption  that $h^0(H)$ or $H^3$ is sufficiently large.

\section{Enlarging the Miyaoka-Mori family}

In the remainder of the paper we will use the notion  of rationally 
connected fibration associated to a covering family of rational curves.
We refer to \cite{Campana1} for the construction (cf. also \cite{KoMiMo1}).

\begin{Theorem}  Suppose $X$ is a 3-dimensional Fano variety with
log-terminal
singularities  of fixed index $n.$ Then for some constant $c$
either $(-K_X)^3 \le c$ or 
there exists a covering  family of rational curves $\{l\}$ on $X$
with the following properties.

1) The degree $l\cdot H \le c$

2) The rationally connected fibration associated to $\{l\}$ has fibers of
dimension $2$.
\end{Theorem}
      
{\bf Proof.}  Suppose $X$ is as above.
By Miyaoka-Mori theorem (cf. \cite{MiyaokaMori}) there is a 
covering family $\{l\}$ of rational curves on $X$ such that
$l\cdot (-K_X) \le 6.$

Consider the associated rationally connected fibration. If its image is a
point,  $(-K_X)^3$ is bounded (cf., e.g. \cite{Fano}). If its image is a curve 
then we are done. So we just need to consider the case when the associated
rationally connected fibration has image of dimension 2. In this case
the family $\{l\}$ is the family of fibers.
If a general $l$ passes through the singularities of $X,$ one can use the
methods of \cite{Fano}, sections 5,6, to construct a new family whose rationally connected
fibration has fibers of dimension at least 2, and degree is bounded in terms
of $n.$ So we can assume that $\{l\}$ doesn't pass through $Sing(X).$

In this case $\{l\}$  is a free family on $X\setminus Sing(X)$  in the terminology
of Nadel (\cite{Nadel}). Consider the associated rationally connected
fibration on the terminal modification $\pi :Y\rightarrow X.$  That is, we
have a Zariski open subset $U$ of $Y$ and a proper morphism 
$\phi : U \rightarrow Z.$ Take a ``general" curve $C$ on $Z$. Consider
the divisor  $D$ which is a Zariski closure of $\phi^{-1}(C).$ By the
definition of $Y$, $D$ is $\Q-$Cartier, $k\cdot D$ is Cartier for some 
$k\in \N.$  Consider an ample $\Q-$Cartier divisor $M$ on $Y$ which is
obtained by subtracting some exceptional divisors of $\pi$ from $\pi^*(-K_X).$

By construction $D\cdot l= 0$ and $(K_X+M)\cdot l =0.$ 
For big enough and divisible enough $m,$ $mM-kD$ is very ample.
Let's denote it by $F$. Then one can apply the Alexeev Minimal Model Program
(cf. \cite{Alexeev1}) to $(Y,K_Y+\frac1{m}|F|)$.

Let us denote the end result of this program by $Y_1$. From the construction there exists a one-dimensional Zariski closed subset $S$ of  $Y_1$ such that the restriction of $(\pi ^{Y}_{Y_1})^{-1}$ to $U=Y_1 \setminus S$ is an isomorphism.

We have several possibilities for $Y_1$. Suppose first that $Y_1$  is a minimal model. Then $K_{Y_1}+\frac1{m}|F|_{Y_1}$ is numerically effective on $Y_1$. By intersecting two general hyperplane sections we can choose a covering family of (not necessarily rational) curves $\{C\}$ such that $C\subset U$ and $C_Y \cdot D >0.$ Then we get a contradiction as follows.
$$0\leq (K_{Y_1}+\frac1{m}|F|_{Y_1})\cdot C = (K_Y+\frac1{m}|F|)\cdot C_Y \leq -\frac1m(kD)\cdot C_Y <0$$
So the $Y_1$ is not a minimal model, and thus we have a Mori fibration $Y_1\rightarrow Z.$

If $\dim Z=0$ then we can immediately bound $(-K_X)^3.$ Indeed, for any positive divisible enough $p$ we have the following inequalities.
$$h^0(-pK_X)\leq h^0(-pK_Y) \leq h^0(-pK_{Y_1})$$
By the Kawamata-Vieweg vanishing and Riemann-Roch, $h^0(-pK_{Y_1})$ grows like $\frac{1}{3!}(-K_{Y_1})^3\cdot p^3$ and $h^0(-pK_X)$ grows like $\frac{1}{3!}(-K_X)^3\cdot p^3.$ By the result of Kawamata \cite{Kawamata} (cf. also \cite{KMMT}) $(-K_{Y_1})^3$ is bounded. Thus $(-K_X)^3$ is bounded.

If $\dim Z =1$ then the fibers are smooth del Pezzo surfaces. Because they are rational and form a bounded family, we can choose a covering family of curves $\{L\}$ on $Y_1$ such that a general $L$ belongs to a general fiber, connects two general points of the fiber, and has bounded intersection with the anti-canonical class. We can then pull it back to $Y$ and down to $X.$ This family would obviously satisfy our requirements. Note also that this would be a coextremal ray on $Y$ in the terminology of Batyrev \cite{Batyrev2}.

If $\dim Z=2$ then we claim that the family of fibers of the fibration is not $\{l\}.$ Suppose it is $\{l\}.$ If at any step of the program we contract a prime divisor intersecting $\{l\}$ then $\{l\}$ would cease to be free and would never become free again. Thus no such contraction ever occured. But this means that
$$(K_{Y_1}+\frac1m F_{Y_1})\cdot l = (K_Y+\frac1m F_Y)\cdot l$$
The left hand side of the above equality is negative, while the right hand side is zero by the definition of $F$ and because $l$ doesn't pass through the singularities of $X.$

So the family of fibers $\{L\}$ is not $\{l\}$. It is a free family on $Y_1$ and $-K_{Y_1}\cdot L=2.$ By pulling it back to $Y$ we get a free family on $Y$ with the same intersection with $-K_Y.$ We can then glue it with $\{l\}$ (see \cite{KoMiMo2}) to get the desired covering family of rational curves. \hfill \checkbox


\section{Completion of the proof}

First of all, let's define some number $M =M(n)$ which will play a crucial role
in the proof of the main theorem.

\begin{Lemma} For every $n$ there exists some constant $M\ge 4320n$ such that for 
any del Pezzo surface with log-terminal singularities of index dividing $n$,  
and any divisor $D\in |-4320nK_S|$, the pair $(S, \frac{1}{M} D)$ is a 
log-terminal pair.   
\end{Lemma}

{\bf Proof.} It follows essentially from the boundedness of del Pezzo 
surfaces of given index. First of all, the multiplicity of any of the 
irreducible components of $D$ can not be too big because $H\cdot D$ is 
bounded for some very ample $H$ on $X$ (which can be actually chosen to be 
some fixed multiple of $(-nK_S)$). So for big enough $M,$ $(S, \frac{1}{M} D)$ 
is a log pair, i.e. $\frac{1}{M} D$ is a boundary.
Moreover, on the minimal resolution
$\pi: Y \rightarrow S$ the exceptional divisors have simple normal crossings. For every such exceptional divisor $E$ there is a covering family of curves $\{L\}$ on $Y$ that intersect $E$, such that $L\cdot (-K_S)$ is bounded. This implies that for
every divisor $D$ as above,  the coefficient of $E$ in $\pi ^*D$ is bounded. So for each $Y$ we have only finitely many possibilities for the linear system that the strict preimage $D'$ of $D$ belongs to. Altogether we have a bounded family of pairs $(Y,\pi ^*D)$. For each pair there exists the minimal number $M(Y,\pi ^*D)$ such that the pair $(S, \frac{1}{M} D)$ is log-terminal. This $M$ can be found by successive blow-ups of $Y$. Clearly, it is Zariski semi-continuous as a function of the pair $(Y,\pi ^*D)$, so it is uniformly bounded.

\hfill \checkbox

{\bf Remark.}  The bound in the above lemma is not explicit. It would be very interesting to find an explicit bound for $M$ in terms of $n$.

In the remainder of this section we consider the Fano threefolds $X$ with log-terminal singularities of index $n$, which are equipped with the family of rational curves $\{l\}$ from the statement of Theorem 3.1. Let us denote the family of fibers of its rationally connected fibration by $\{S\}$. This is a LINEAR system of divisors (cf. \cite{Fano}). 

We have the following easy proposition.
\begin{Proposition} In the above situation and notations, there exist constants $c_1$ and $c_2$ depending only on the degree $H\cdot l$ such that for a general $S$

$(H_{|_S})^2 \le c_1,$ and $h^0(H_{|_S})\le c_2.$
\end{Proposition}

{\bf Proof.} Possibly gluing $\{l\}$ with itself once, we may assume that two general points of $S$ are connected by some curve $l$. Suppose $l\cdot H= m.$ 
By Lemma 2.2 of \cite{Fano}, $(H_{|_S})^2 \le m^2.$ To bound $h^0(H_{|_S}),$ consider $(m+1)$ curves $l$ and $(m+1)$ points on each of them. If a section of 
$H^0(H_{|_S})$ vanishes at all these points, it has to vanish at all these curves. This implies that it is identically zero on $S$, by intersecting with a general $l,$ because $l\cdot l \ge 1.$ So $h^0(H_{|_S})\le (m+1)^2$.

\begin{Theorem} In the above situation and notations, suppose that the base locus of $|S_X|$ contains a curve. Then $H^3$ is bounded.
\end{Theorem}

{\bf Proof.} The divisor $H$ is very ample on $X.$ Consider the family of (probably not rational) curves $\{C\}= H_{|_S}.$ By the above proposition, $C\cdot H$ is bounded. Because $X$ has log-terminal singularities of index  $n$, and the base locus of $|S_X|$ contains a curve, we can find a point $x_0$ in this base locus, with multiplicity of the local ring bounded by $2n$. Then the curves $C$ passing through this point are not contained in any closed subvariety of $X.$ So $H^3$ is bounded (cf., e.g. \cite{Fano}, Lemma 2.2.) \hfill \checkbox

{\bf Remark.}  If $X$ is $\Q-$factorial, and $\rho(X)=1,$ the conditions of the
above theorem are always satisfied. So in this case one can use the above theorem instead of the argument of \cite{Fano}, section 4, to get a bound for $(-K_X)^3$, which is polynomial in $n.$

For the remainder of this section we will assume that the base locus of $S_X$ contains no curves.  If $X$ is $\Q-$factorial, this implies that it is empty, because it is equal to the intersection of two general elements of the one-dimensional linear system $S_X$. In general, it could consist of finitely many points.

By the Log Minimal Model Program (cf., e.g. \cite{FlipsAbundance}, 6.16; also \cite{Shokurov1}) there exists a small partial resolution of singularities $\pi: Y\rightarrow X$ such that $Y$ is $\Q-$factorial, and $S$ is $\pi -$nef. (If $X$ is $\Q-$factorial then $Y=X$).
Because $S_Y$ is $\pi -$nef, the base locus of $S_Y$ contains no exceptional curves of $\pi ,$ so it is empty. Therefore, $|S_Y|$ is a free linear system,
two general $S_Y$ do not intersect. The variety $Y$ has log-terminal singularities of index $n$. The anti-canonical class $-K_Y= \pi^*(-K_X).$
So by trivial adjunction a general $S_Y$ is a del Pezzo surface with log-terminal singularities of index dividing $n.$

By the Cone Theorem on $Y$, we have at least one of the following.
 
\begin{enumerate}
\item One of the extremal curves of $\pi,$  say $L$, intersects $S_Y$.

\item There is an extremal contraction on $Y$ whose fibers intersect $S_Y.$

\end{enumerate}

In the second case we could have a fibration, a divisorial contraction, or a small contraction. In the fibration case, we can get a covering family of curves on the fibers with bounded intersection with $-K_Y,$ whose rationally connected fibration is different from $|S|$. So we can glue together this family and $\{l\}$ to get a new family, which connects two general points of $X$ and has bounded intersection with $H$. This implies boundedness of $H^3,$ so we are left with the following three possibilities.

\begin{enumerate} 

\item There is a divisorial contraction on $Y$ whose fibers intersect $S$.

\item There is an exceptional curve $L$ on $Y$, such that $L\cdot K_Y<0,$ and
$L$ intersects $S$.  

\item One of the extremal curves of $\pi,$  say $L$, intersects $S_Y$. In this case $L\cdot K_Y =0$.

\end{enumerate}

The first case is the simplest. It is treated in the following proposition.

\begin{Proposition} Suppose $X$ and $Y$ are as above, and we have a divisorial contraction on $Y$. Then $h^0(H)$ is bounded.
\end{Proposition}

{\bf Proof.}
Suppose $E$ is the exceptional divisor of the contraction. As in \cite{Fano}, section 6, we can find a covering family $\{L\}$ on $Y$ such that $-3\le K_{Y_1}\cdot L \le 0$ and $-E\cdot L \le 3.$ (The variety $Y$ in our case has log-terminal singularities rather than terminal, but the Subadjunction still works).

This implies that $L\cdot H$  is bounded. Suppose $L\cdot H=d.$ By Proposition 4.1, $h^0(H_{|_S})$ is bounded. Suppose it is equal to $e.$ Also fix a constant $M$ from Lemma 4.1. Then  $h^0(H) - h^0(H_Y-(3M+d+1)S_Y) \le (3M+d+1)e$. If $h^0(H)$ is big enough, then there exists some $D\in |H_Y-(3M+d+1)S_Y)|.$ Let us write 
$D= a\cdot E + F,$ where $E\notin SuppF.$
By Lemma 4.1 $(S, \frac{1}{M}D_{|_S})$ is a log pair. This means that $\frac{1}{M}D_{|_S}$ is a boundary. Therefore $a \le M.$ The following chain of inequalities provides a contradiction.
$$d= L\cdot H= (3M+d+1)L\cdot S + L\cdot D \ge (3M+d+1)L\cdot S + aL\cdot E  \ge$$
$$\ge (3M+d+1)L\cdot S -3M \ge (d+1) $$ 
\hfill \checkbox

Note that in this case we did not use the full strength of Lemma 4.1. In particular in this case we could actually choose $M$ to be polynomial in $n$.
We will however need the full strength of the lemma for the remaining two cases.
We will treat them together in the following proposition.

\begin{Proposition} Suppose $X$ and $Y$ are as above, and we have a small 
contraction on $Y$ with some fiber $L$ intersecting $S$. 
Then $h^0(H)$ is bounded.
\end{Proposition}

{\bf Proof.}  Suppose $L$ is as above. Recall again a number $M$ from Lemma 4.1. By Proposition 4.1, if $h^0(H_X) = h^0(H_Y)$ is big enough then $h^0(H_Y-(3M)S_Y) >0. $ Take some 
$F \in |H_Y-(3M)S_Y|$ and consider the divisor $D$ which consists of all the components of $F$ containing $L.$ 

By the definition of $M$ the log pair $(Y, \frac{1}{M}D )$ is log-terminal in the neighborhood of the general point of $L.$  Consider the log-terminal modification $(Y_1, B_1)$ of $(Y, B)$, where $B=\frac{1}{M}D $ (cf. \cite{FlipsAbundance}, 6.16, d=1). The log adjunction formula for the morphism $\pi^{Y_1}_{Y}$ is the following.
$$K_{Y_1} +B_1 = (\pi^{Y_1}_{Y})^* (K_Y+B) +\sum \alpha_i A_i,$$
where $A_i$ are exceptional divisors of  $\pi^{Y_1}_{Y}$ and all $\alpha_i\le 0$
because $K_{Y_1} +B_1$  is relatively nef.

Because $(Y, \frac{1}{M}D )$ is log-terminal in the neighborhood of the 
general point of $L,$ none of the divisors $A_i$ lies above $L.$ So we can 
consider the curve $L_1$ on $Y_1$ which is the closure of the pullback of 
the general point of $L.$  Obviously, $L_1$ is an extremal curve on $Y_1,$ and 
$L_1\cdot ( K_{Y_1} +B_1) \le L\cdot (K_Y+B).$ By Kawamata's theorem on the 
length of an extremal curve (cf. \cite{Ka2}), $-2\le L_1\cdot ( K_{Y_1} +B_1).$

So we have the following chain of inequalities.

$$-2\le L_1\cdot ( K_{Y_1} +B_1) \le L\cdot (K_Y+B) \le L\cdot 
(K_Y +\frac{1}{M}F) =$$
$$=L\cdot (K_Y +\frac{1}{M}(H_Y -3M\cdot S)) = L\cdot K_Y \cdot (1-\frac{4320n}{M} ) -3L\cdot S \le -3$$ 
The above contradiction completes the proof. \hfill \checkbox


\end{document}